\documentstyle[amscd,amssymb,12pt]{amsart}
\input xy
\xyoption{all}

\CompileMatrices                                                             
\newdir{((}{{\lhook\kern-.1em}}                                              
\newdir{))}{{\rhook\kern-.1em}}                                              
\newdir{ >}{{}*!/-5pt/@{>}}                                                  

\newcommand{\la}{\langle}
\newcommand{\ra}{\rangle}

\newcommand{\ve}{\varepsilon}

\newtheorem{theorem}{\bf Theorem}[section]
\newtheorem{lemma}[theorem]{\bf Lemma}

\newtheorem{corollary}[theorem]{\bf Corollary}

\newcommand{\be}{\begin{equation}}
\newcommand{\ee}{\end{equation}}        

\newfont{\bfc}{cmbsy10 scaled 1200}  
\newfont{\dr}{msbm10 scaled \magstep1}  
\newfont{\sdr}{msbm8}  
\newfont{\gl}{eufm10 scaled \magstep1}  

\DeclareFontFamily{OT1}{rsfs}{}
\DeclareFontShape{OT1}{rsfs}{n}{it}{<->rsfs10}{}
\DeclareMathAlphabet{\curly}{OT1}{rsfs}{n}{it}

\newcommand{\CC}{{\Bbb C}}
\newcommand{\DD}{{\Bbb D}}

\newcommand{\QQ}{{\Bbb Q}}
\newcommand{\RR}{{\Bbb R}}
\newcommand{\TT}{{\Bbb T}}
\newcommand{\ZZ}{{\Bbb Z}}

\newcommand{\glie}{{\frak g}}

\newcommand{\klie}{{\frak k}}

\newcommand{\tlie}{{\frak t}}

\newcommand{\fX}{{\curly X}}

\newcommand{\Ad}{\operatorname{Ad}}

\newcommand{\Ext}{\operatorname{Ext}}

\newcommand{\Ker}{\operatorname{Ker}}
\newcommand{\Lie}{\operatorname{Lie}}
\newcommand{\Map}{\operatorname{Map}}
\newcommand{\Mon}{\operatorname{Mon}}

\newcommand{\Tor}{\operatorname{Tor}}

\renewcommand{\exp}{\operatorname{exp}}
\newcommand{\id}{\operatorname{id}}
\renewcommand{\Im}{\operatorname{Im}}

\newcommand{\AAA}{{\curly A}}

\newcommand{\imag}{i}
\newcommand{\qu}{/\kern-.7ex/}
\newcommand{\exh}{\to\kern-1.8ex\to}

\title[Lifts of smooth group actions to line bundles]
{Lifts of smooth group actions to line bundles}
\author{Ignasi Mundet i Riera}
\address{Centre de Math{\'e}matiques \\
{\'E}cole Polytechnique \\
Palaiseau, France \\
and
Departamento de Matem{\'a}ticas \\
Universidad Aut{\'o}noma de Madrid \\ 
Madrid, Spain}
\date{10--1--2000}
\email{ignasi@@math.polytechnique.fr}

\begin{document}
\maketitle
\begin{abstract}
Let $X$ be a compact manifold with a smooth action of a compact
connected Lie group $G$. Let $L\to X$ be a complex line bundle. Using
the Cartan complex for equivariant cohomology, we give 
a new proof of a theorem of Hattori and Yoshida which
says that the action of $G$ lifts to $L$ if and only if the first
Chern class $c_1(L)$ of $L$ can be lifted to an integral equivariant  
cohomology class in $H^2_G(X;\ZZ)$, and that the different lifts
of the action are classified by the lifts of $c_1(L)$ to $H^2_G(X;\ZZ)$.
As a corollary of our method of proof, we prove that, if the action is
Hamiltonian and $\nabla$ is a connection on $L$ which is unitary for
some metric on $L$ and whose curvature is $G$-invariant, 
then there is a lift of the action to a certain
power $L^d$ (where $d$ is independent of $L$)
which leaves fixed the induced metric on $L^d$ and 
the connection $\nabla^{\otimes d}$. This
generalises to symplectic geometry a well known result in Geometric
Invariant Theory.
\end{abstract}

\section{Introduction and statement of the results}
Let $X$ be a connected smooth compact manifold with a smooth 
left action of a compact connected Lie group $G$. 
Our aim is to study liftings of the action of $G$ to
complex line bundles $L\to X$. Of course, it is not always possible to
find such a lift. For example, if $x\in X$ has trivial stabiliser,
then the restriction $L|_{Gx}$ has to be topologically trivial.

The problem which we consider is very natural and has been already studied 
by several people. The general question on lifting of smooth actions
to principal bundles was considered for example by R. Palais and T. Stewart
\cite{PS, S}. The more concrete problem of lifting actions to complex
line bundles was studied for example by B. Kostant in \cite{Ko}, where
he proved that if $G$ is simply connected, $X$ is symplectic and the
action of $G$ is Hamiltonian, then there is always some lift of the action to
$L$ (see Theorem 4.5.1 in \cite{Ko}), and by A. Hattori and T. Yoshida
\cite{HY} and R. Lashof \cite{L}.

Let $EG\to BG$ be the universal $G$-principal bundle.
Fix a point $x_0\in BG$, and denote by
$\iota:X\to X_G$ the inclusion of the fibre over $x_0\in BG$ 
of the Borel construction $X_G=EG\times_G X\to BG$.

Our first result is a new proof of the following theorem, which was
proved in \cite{HY}. Our method of proof is however different from theirs.

\begin{theorem}
Let $L\to X$ be a line bundle. 
The action of $G$ on $X$ lifts to a linear action on $L$ if and only
if $$c_1(L)\in \iota^* H^2_G(X;\ZZ).$$
Furthermore, if $c_1(L)\in\iota^*H^2_G(X;\ZZ)$, then the different lifts of the action are classified by $\iota^{-1}(c_1(L))$.
\label{lifting}
\end{theorem}

Fix from now on a metric on $L$. All the connections we will take on
$L$ will be assumed to be unitary with respect to this metric, and all
the actions of $G$ on $L$ will keep the metric fixed.

Let $\glie=\Lie(G)$.
To describe lifts of the action of $G$ and invariant connections on
$L$ we will use the Cartan model for real equivariant cohomology (see
section \ref{cohequi} for the necessary definitions).
Suppose that $c_1(L)=\iota^* l$ for some $l\in H^2_G(X;\ZZ)$, and let
$\alpha-\mu$ be a closed element in the Cartan complex
$\Omega^*_G(X;\imag\RR)$ representing the class 
$-2\pi\imag l\in H^2_G(X;\imag\RR)$,
where $\alpha\in\Omega^2(X;\imag\RR)^G$ and
$\mu\in\Omega^0(X;\imag\glie^*)^G$. 
Any connection $\nabla$ on $L$ whose curvature is $\alpha$ can be
combined with $\mu$ to obtain an infinitesimal lift of the action 
(see section \ref{inflift} for definitions and Theorem
\ref{tinflift}), and we will study whether there is a connection $\nabla$
which defines an infinitesimal lift that exponentiates to an action of
$G$. Let $\TT^{\alpha}$ be the gauge equivalence classes of connections on
$L$ whose curvature is $\alpha$, and let $a_1:H_1(G;\ZZ)\to H_1(X;\ZZ)$ be the
map induced by the map 
$a_1(x):G\ni g\mapsto gx\in X$ for any $x$ (since $X$ is 
connected, the map in homology is independent of $x$).

\begin{theorem} The set $\TT^G_{\alpha}(\mu)$ of gauge 
equivalence classes of connections 
with curvature $\alpha$ which, combined with $\mu$, define an 
infinitesimal lift which
exponentiates to an action of $G$, is 
a subtorus of $\TT^{\alpha}$ of dimension $b_1(X)-\dim(\Im
a_1\otimes_{\ZZ}\RR)$, where $b_1(X)=\dim H^1(X;\RR)$. 
In particular, if $X$ is symplectic and the action of $G$ is
Hamiltonian, then $\TT^G_{\alpha}(\mu)=\TT^{\alpha}$.
\label{descripcio}
\end{theorem}

In fact, we get in this way all the possible lifts of the action. This is due 
to the following reasons: (1) given any lift, there is some invariant 
connection (just take any connection and average); (2) using the construction
of differential forms in the Cartan complex representing the first equivariant
Chern class $c_1^G(L)$ of $L$ (see \cite{BV}) we get a closed element 
$\alpha-\mu$ representing the class $-2\pi\imag c_1^G(L)$; and (3), the group 
$G$ is connected, so a lift of the action is uniquely determined by the 
infinitesimal action on $L$.

The results and techniques in this paper have interesting consequences
in the case in which $X$ is symplectic and the action of $G$ is Hamiltonian.
The following corollary is an analogue of a well known result in Geometric
Invariant Theory (see Corollary 1.6 in \cite{MFK}).

\begin{corollary}
Let $G$ act on a symplectic manifold $(X,\omega)$ in a Hamiltonian fashion.
Then there exists an integer $d\geq 1$ with the following
property. For any line bundle $L\to X$ with a connection $\nabla$
whose curvature is $G$ invariant,
there is a lift of the action of $G$ to $L^d$ such that the induced
connection $\nabla^{\otimes d}$ on $L^d$ is $G$-invariant.
\label{simplectic0}
\end{corollary}

\begin{corollary}
Let $G$ act on a symplectic manifold $(X,\omega)$ in a Hamiltonian fashion.
Let $\ve>0$.There exists a symplectic form $\omega'$ such that

(1) $\omega'$ is preserved by $G$,

(2) $|\omega-\omega'|_{C^0}<\ve$,

(3) the action of $G$ on $(X,\omega')$ is Hamiltonian,

(4) there is a natural number $k>0$, a Hermitian line bundle $L\to X$ 
with a unitary connection $\nabla$ with curvature $\imag k\omega'$ and a
linear action of $G$ lifting the action on $X$ and preserving $\nabla$.
\label{simplectic}
\end{corollary}

This paper is organized as follows: in section 2 we recall some facts
on equivariant cohomology which we will need; in section 3 we state the
relation between infinitesimal lifts and 2-forms in the Cartan
complex; in section 4 we define and prove some key properties of the 
monodromy map $M_{\gamma}$, which
measures how far an infinitesimal lift is to exponentiate to an action
of $G$; in section 5 we study how to chose connections which provide
liftings of the action of $G$; finally, in section 6 we give the proofs of
Theorems \ref{lifting} and \ref{descripcio}, and Corollaries
\ref{simplectic0} and \ref{simplectic}. 

\noindent{\bf Acknowledgements.}
I thank M. Vergne for having pointed out to me the paper of Kostant \cite{Ko}.
I also thank O. Garc{\'\i}a--Prada for some comments on the paper.

\section{Equivariant cohomology}
\label{cohequi}

In this section we recall some basic facts on equivariant
cohomology. For more information the reader is refered to 
\cite{AB, BGV, GS}.
 
Let $\pi:X_G=EG\times_G X\to BG$ be the Borel construction of $X$.
The equivariant cohomology of $X$ is by definition the singular
cohomology of $X_G$, and is denoted, for any ring $R$, as
$$H^*_G(X;R):=H^*(X_G;R).$$

Let $\glie$ be the Lie algebra of $G$.
Let $\Omega^*(X)$ be the complex of differential forms on $X$.
Denote by $\fX:\glie\to\Gamma(TX)$ the map which assigns to any
$s\in\glie$ the vector field on $X$ generated by the infinitesimal
action of $s$. So, if $f\in\Omega^0(X)$, $x\in X$ and $s\in\glie$,
then $\fX(s)(f)(x)=\lim_{\ve\to 0}\ve^{-1}(f(e^{\ve s}x)-f(x))$.
Consider the complex
$$\Omega^*_G(X)=(\Omega^*(X)\otimes\RR[\glie])^G$$
(as usual the supscript $^G$ means $G$-invariant elements) with the
grading obtained from the usual grading in $\Omega^*(X)$ and twice the
grading in $\RR[\glie]$ given by the degree, together
with the differential $d_\glie$ defined by 
$$d_\glie(\eta)(s)=d(\eta(s))+\iota_{\fX(s)}\eta(s),$$
where $\eta\in\Omega^*_G(X)$, $s\in\glie$, and
$\iota_v:\Omega^*(X)\to\Omega^{*-1}(X)$ is the contraction map.
One can check that $d_{\glie}^2=0$, and the complex 
$(\Omega^*_G(X),d_{\glie})$ is called the Cartan complex.
It will be our main tool in this paper.
Note that in this paper we consider sometimes elements of
$\Omega^*_G(X;\imag\RR)=\imag\Omega^*_G(X)$. 
 
The following classical theorem (which is a generalisation of de
Rham's theorem) is proved for example as Theorem 2.5.1 in \cite{GS}.

\begin{theorem}
There is a natural isomorphism 
$H^*_G(X;\RR)\simeq H^*(\Omega^*_G(X),d_{\glie})$.
\end{theorem}

We will need an explicit description of the isomorphism given by the
theorem in degree 2,
as given by the following lemma. The proof can easily be deduced
from the proof of Theorem 2.5.1 in \cite{GS}.

\begin{lemma}
Let $\alpha-\mu\in\Omega^2_G(X)$ satisfy $d_{\glie}(\alpha-\mu)=0$. 
Let $f:\Sigma\to X_G$ be a
continuous map. Let $P_f:=(\pi\circ f)^*EG\to \Sigma$, and let $\phi_f$
be the induced section of $\pi_f:X_f=P_f\times_G X=(\pi\circ
f)^*X_G\to\Sigma$. Suppose that both $P_f$ and $\phi$ are smooth, and
let $A$ be a connection on $P_f$. We then get a projection
$\pi_A:TX_f\to \Ker d\pi_f$,
which allows to pullback $\alpha$ to a form
$\pi_A^*\alpha\in\Omega^2(X_f)$ vanishing on the tangent vectors which
are horizontal with respect to $A$. Then
$$\la [\Sigma],f^*[\alpha-\mu]\ra =
\int_{\Sigma}\phi_f^*(\pi_A^*\alpha)-\la \mu,F_A\ra,$$
where $F_A\in\Omega^2(P_f\times_{\Ad}\glie)$ is the curvature of $A$.
\label{mapequi}
\end{lemma}

Observe that a constant and equivariant map $\mu:X\to\glie^*$
represents an integral element $[\mu]\in H^2_G(X;\ZZ)$ if and only if
for any $u\in\glie$ such that $\exp(u)=1$ we have $\la\mu,u\ra\in\ZZ$.

Finally, the map $\iota^*:H^2_G(X;\RR)\to H^2(X;\RR)$ can be written
as follows using the Cartan and de Rham complexes: if $\alpha-\mu$ is
a closed element of $\Omega^2_G(X)$, then 
\begin{equation}
\iota^*[\alpha-\mu]=[\alpha]\in H^2(X;\RR)
\label{mapiota}
\end{equation}
(of course, since $d_{\glie}(\alpha-\mu)=0$ we have $d\alpha=0$).

\section{Infinitesimal lift of the action}
\label{inflift}
Let $\nabla$ be a connection on $L$ whose
curvature $\alpha=\nabla^2\in\Omega^2(X;\imag\RR)$ is $G$-invariant.
We will call an infinitesimal lift of the action of $G$ on $X$ to
an action on $L$ any map $\widetilde{\fX}:\glie\to\Gamma(TL)$ which
satisfies\footnote
{Since the action of $G$ is on the left, for any
$s,s'\in\glie$ we have
$\fX([s,s'])=-[\fX(s),\fX(s')].$}
$\widetilde{\fX}([s,s'])=-[\widetilde{\fX}(s),\widetilde{\fX}(s')]$
for any $s,s'\in\glie$ and such that
$$d\pi\circ\widetilde{\fX}=\fX.$$

Let $u=2\pi\imag\in\imag\RR=\Lie(S^1)$. 
Let $U_L\in\Gamma(\Ker d\pi)$ be the vertical tangent
field generated by the infinitesimal action of $s$ on $L$ given by
fibrewise multiplication.

Let $\fX^{\nabla}:\glie\to\Gamma(TL)$ be the map which assigns to $s\in\glie$
the horizontal lift of $\fX(s)$ obtained using $\nabla$. The following
is well known (see for example Section 3 in \cite{Ko}).

\begin{theorem}
Let $\mu:X\to\imag\glie^*$ be a map which satisfies
$d\mu(s)=\iota_{\fX(s)}\alpha$ for any $s\in\glie$.
Then $\widetilde{\fX}^{\nabla,\mu}
:=\fX^{\nabla}+\imag\mu U_L$ is an infinitesimal lift
of the action of $G$ on $X$.

Conversely, for any infinitesimal lift $\widetilde{\fX}$ which leaves
$\nabla$ invariant, the function $\mu:X\to\imag\glie^*$ defined by
$\widetilde{\fX}:=\fX^{\nabla}+\imag\mu U_L$ satisfies 
$d\mu(s)=\iota_{\fX(s)}\alpha$ for any $s\in\glie$.
\label{tinflift}
\end{theorem}

Note that the condition 
$d\mu(s)=\iota_{\fX(s)}\alpha\text{ for any }s\in\glie$
is equivalent to asking $$d_{\glie}(\alpha-\mu)=0.$$

\section{The monodromy map $M_{\gamma}$}

Let $(L,\nabla)\to X$ 
be as in the preceeding section, let $\alpha=\nabla^2$ be the
curvature of $\nabla$, and let $\mu:X\to\imag\glie^*$ be 
a map which satisfies $d\mu(s)=\iota_{\fX(s)}\alpha$ for any $s\in\glie$.
Let $\widetilde{\fX}=\widetilde{\fX}^{\nabla,\mu}
:\glie\to\Gamma(TL)$ be the corresponding
infinitesimal lift.

Let $\gamma:S^1\to G$ be any representation.
As before, let $u=2\pi\imag\in\imag\RR=\Lie(S^1)$.
For any $x\in X$ and $y\in L_x$, 
let $\nu_x:[0,1]\to L$ be the integral line of the
vector field $\widetilde{\fX}(\gamma_*(u))$ with initial value $\nu_x(0)=y$.
We then have $\nu_x(1)\in L_x$, and there is a unique
$M_{\gamma}(x)=M_{\gamma}^{\nabla,\mu}(x)\in S^1$ 
(independent of $y$) such that $\nu_x(1)=M_{\gamma}(x)\nu_x(0)$.

The map $M_{\gamma}:X\to S^1$ which we have defined measures the
extent to which the infinitesimal action given by $(\nabla,\mu)$
exponentiates to an action of the group. 
We will see below (viz. Lemma \ref{sasufi}) that the
condition $M_{\gamma}=1$ for any $\gamma$ is enough to ensure that the
infinitesimal action exponentiates.

Let $x\in X$ be any point, and let 
$\Mon^{\nabla}(S^1x)\in S^1$ be the monodromy of parallel transport using
$\nabla$ along the path
$[0,1]\ni t\mapsto \gamma(e^{2\pi\imag t})x$.
The following formula for $M$ can be easily proved using coordinates
in a neigbourhood of $S^1x$ (see also Theorem 2.10.1 in \cite{Ko}):
\begin{equation}
M_{\gamma}^{\nabla,\mu}(x)=\Mon^{\nabla}(S^1x)
\exp(-2\pi\la\mu(x),\gamma_*(u)\ra).
\label{moncon}
\end{equation}
An easy consequence of this formula is that $M$ is gauge invariant,
i.e., for any gauge transformation $g:L\to L$,
\begin{equation}
M_{\gamma}^{g^*\nabla,\mu}(x)=M_{\gamma}^{\nabla,\mu}(x).
\label{gaugeinvariant}
\end{equation}
Let $\eta\in\Omega^1(X;\imag\RR)$, so that $\nabla+\eta$ is another
connection on $L$. Then we also deduce from (\ref{moncon}) that
\begin{equation}
M^{\nabla+\eta,\mu}_{\gamma}(x)=
M^{\nabla,\mu}_{\gamma}(x)\exp\left(-\int_{S^1}\gamma_x^*\eta\right),
\label{addicio}
\end{equation}
where $\gamma_x:S^1\to X$ maps $\theta$ to $\gamma(\theta)x\in X$.

\subsection{Cohomological interpretation of $M$}
In this subsection we assume $G=S^1$ and $\gamma=\id$, and we
denote $M=M_{\gamma}$. 
We identify $\imag\RR$ with $(\imag\RR)^*$ by assigning to
$\alpha\in\imag\RR$ the map $\imag\RR\ni a\mapsto
\la\alpha,a\ra=-\alpha a/2\pi$.
In particular, $\mu\in\Omega^0(X;\RR)^{S^1}$.
Let us suppose that the action of $S^1$ on $X$
is generically free, i.e., the isotropy group is trivial for generic
$x\in X$ (if this is not true, then either the action of $S^1$ on $X$
is trivial, and
the results in this section are obvious, or there is a biggest common
stabiliser $\{1\}\neq Z\subset S^1$, which acts freely on $X$, and we 
replace $X$ by $X/Z$ and $S^1$ by $S^1/Z$).

\begin{lemma}
Let $x\in X$ be a point with trivial stabiliser, and let $O_x:S^1\to
X$ be the map $O_x(\theta)=\theta x$. Assume that $O_x(S^1)$ is
homologous to zero. Then, if $\frac{\imag}{2\pi}[\alpha-\mu]\in 
H^2_{S^1}(X;\ZZ)$, we have $M(x)=1$. 
\label{homognul}
\end{lemma}
\begin{pf}
Let $\Sigma_0$ be a compact surface with a fixed isomorphism 
$\partial\Sigma_0\simeq S^1$ and compatible orientation, 
and let $b_0:\Sigma_0\to X$ be a map
such that $b_0|_{\partial \Sigma_0}=O_x$. Note that $\iota\circ
b_0(\partial \Sigma_0)=\iota\circ O_x(S^1)$ is contained in
$ES^1\times_{S^1}(S^1x)\subset X_{S^1}$. Since the action of $S^1$ on
$S^1x\subset X$ is free,  $ES^1\times_{S^1}(S^1x)$ is
contractible. So, denoting the unit disk by $\DD$,
we may take a map $B_1:\DD\to ES^1\times_{S^1}(S^1x)$ such
that $c\circ B_1|_{S^1}=\iota\circ b_0|_{\partial\Sigma_0}$, where
$c(\theta)=\theta^{-1}$ for $\theta\in S^1$ (such a
$B_1$ is unique up to homotopy rel $S^1$). 
Let us patch together the maps $B_0:=\iota\circ b_0$ and $B_1$
to get a map $B:\Sigma:=\Sigma_0\cup_{S^1}(-\DD)\to X_{S^1}$, where
the minus sign refers to inversed orientation, so that we take the
isomorphism $c$ to identify $\partial(-\DD)$ with
$S^1$. We claim that 
\begin{equation}
M(x)=\exp\la[\Sigma],B^*[-\alpha+\mu]\ra.
\label{forM}
\end{equation}
Clearly, once we check the claim, the lemma is proved.
To prove (\ref{forM}), we will apply Lemma \ref{mapequi} to the map
$B$. Let $\pi:X_{S^1}\to BS^1$ denote the projection. Then
we have $\pi\circ B_0=\{x_0\}$, so that $P_B|_{\Sigma_0}$ (we use here
the notation of Lemma \ref{mapequi}) is the
trivial bundle. Let us fix a trivialisation
$P_B|_{\Sigma_0}\simeq \Sigma_0\times S^1$
of it.  
Using the induced trivialisation $X_B|_{\Sigma_0}=\Sigma_0\times X$,
the restriction $\phi_B|_{\Sigma_0}$ is given by $(\id,b_0)$. 
Since $\DD$ is contractible, $P_B$ is obtained by patching
the trivial bundles over $\Sigma_0$ and $\DD$ through a gluing map
$\rho:\partial\Sigma_0=S^1\to S^1$.
The map $b_0|_{\partial\Sigma_0}$ has winding number $1$, so the section 
$\phi_B|_{\Sigma_0}$ will only glue with a section of the trivial
bundle $(-\DD)\times (S^1x)\subset (-\DD)\times X$ if the gluing map
$\rho$ is the identity $\rho(\theta)=\theta$. Hence, the bundle $P_B$
must have degree $-1$. 

Take now a connection $A$ on $P_B$ which
coincides over $\Sigma_0$ with the flat
connection induced by the chosen trivialisation of $P_B|_{\Sigma_0}$. 
Then the curvature of $A$ is supported in
$\DD$. By Lemma \ref{mapequi}, the RHS of (\ref{forM}) is equal to
$$\int_{\Sigma}-\phi_B^*(\pi_A^*\alpha)+\la \mu,F_A\ra.$$
Now, by our choice of $A$, $\int_{\Sigma}\la \mu,F_A\ra=
\int_{\DD}\la \mu,F_A\ra$. 
On the other hand, the image of the restriction
$\phi_B|_{\DD}$ is contained in $P_B\times_{S^1}(S^1x)\subset X_B$,
and since $\mu$ is $S^1$ equivariant, we may write 
$\int_{\DD}\la \mu,F_A\ra=-\frac{1}{2\pi}\mu(x)\int_{\DD}F_A$. Finally, by
Chern--Weil $\int_{\DD}F_A=2\pi\imag\deg(P_B)=2\pi\imag$.

Using again that $\phi_B(\DD)\subset P_B\times_{S^1}(S^1x)$, we deduce
that $\phi_B^*(\pi_A^*\alpha)$ vanishes on $\DD$. Indeed, the vertical
part of any tangent vector to $P_B\times_{S^1}(S^1x)$ lies in $\Ker
d\pi_B|_{P_B\times_{S^1}(S^1x)}$, which is a real line bundle, so for
any $y\in\DD$ and $u,v\in T_y\DD$, $\pi_A(d\phi_B(u))$ and
$\pi_A(d\phi_B(v))$ are linearly dependent, and hence
$\alpha(\pi_A(d\phi_B(u)),\pi_A(d\phi_B(v)))=0$.
So $\int_{\Sigma}-\phi_B^*(\pi_A^*\alpha)=
\int_{\Sigma_0}-\phi_B^*(\pi_A^*\alpha)$. And this is equal to 
$\int_{\Sigma_0}-b_0^*\alpha$. 
Taking  into account that $\alpha$ is the
curvature of a connection $\nabla$ on the line bundle $L\to X$, 
one can check that
$$\exp\left(\int_{\Sigma_0}-b_0^*\alpha\right)=
\Mon^{\nabla}(S^1x)$$
(see Theorem 1.8.1 in \cite{Ko}).
Now, using (\ref{moncon}) and the preceeding computations we deduce
\begin{align*}
M(x) &=
\exp\left(-\left(\int_{\Sigma}\phi_B^*(\pi_A^*\alpha)\right)-\imag\mu(x)\right)
\\
&=\exp\left(\int_{\Sigma}-\phi_B^*(\pi_A^*\alpha)+\la\mu(x),F_A\ra\right)
\quad\quad\text{ (since $\int F_A=2\pi\imag$)}\\
&=\exp\la[\Sigma],B^*[-\alpha+\mu]\ra.
\end{align*}
\end{pf}

\begin{corollary}
Assume that for some $x\in X$, $O_x(S^1)$ is
homologous to zero. Then, if $\frac{\imag}{2\pi}[\alpha-\mu]\in 
H^2_{S^1}(X;\ZZ)$, we have $M=1$.
\end{corollary}
\begin{pf}
Indeed, the condition of $O_x(S^1)$ being homologous to zero is
independent of $x$, the
function $M$ is continuous, and by assumption the set of $x\in
X$ with trivial stabiliser is dense. 
\end{pf}

When the orbit $S^1x$ is not homologous to zero, the map $M(x)$ will
depend on the connection $\nabla$ (and not only on $\alpha$ and
$\mu$), as we will see below.

\begin{lemma}
Let $x,x'\in X$ be two points with trivial stabiliser. Then
$M(x)=M(x')$. 
\end{lemma}
\begin{pf}
This can be proved either with local coordinates or using the same
technique as above. We sketch the second strategy. For that, let 
$\rho:[0,1]\to X$ be a path such that $\rho(0)=x$ and $\rho(1)=x'$, let
$\Sigma_0=[0,1]\times S^1$ and let $b_0:\Sigma_0\ni (t,\theta)\mapsto 
\theta\rho(t)$. Glue two disks $\DD_0$ and $\DD_1$ to the boundary of
$\Sigma_0$ with suitable orientations to get a closed oriented 
surface $\Sigma$,
and extend the map $\iota\circ b_0$ to a map $B:\Sigma\to X_{S^1}$ 
just as in the preceeding lemma (i.e., so that
the image of $\DD_0$ is contained in $ES^1\times_{S^1}(S^1x)$
and that of $\DD_1$ in $ES^1\times_{S^1}(S^1x')$). As before one
can check that 
$$M(x)-M(x')=\exp\la[\Sigma],B^*[-\alpha+\mu]\ra.$$
Now, however, the map $B$ is homotopic to the trivial map,
and from this the result follows. 
\end{pf}

\begin{corollary}
The map $M:X\to\RR$ is constant.
\label{Mconstant}
\end{corollary}

For the last lemma of this section, we return to the general
situation, in which $G$ is any compact connected Lie group.

\begin{lemma}
Let $\gamma:S^1\to G$ be a morphism, and let $g\in G$. We then have
$$M_{\gamma}=M_{g\gamma g^{-1}}.$$
\label{conjugats}
\end{lemma}
\begin{pf}
Let $\rho:S^1\to G$ be a smooth map such that $\rho(1)=1$ and
$\rho(-1)=g$. Consider on $X\times S^1$ the action of $S^1$ given by
$$\theta(x,\alpha)=\rho(\alpha)\gamma(\theta)\rho(\alpha)^{-1}x
\text{ for $\theta\in S^1$ and $(x,\alpha)\in X\times S^1$.}$$
Let $\pi_1:X\times S^1\to X$ be the projection, and take on $X\times
S^1$ the bundle $\pi_1^*L$ with the connection
$\nabla_{M\times S^1}=\pi_1^*\nabla$. Finally, let $\mu_{M\times
  S^1}(x,\alpha)=\la\mu(x),\Ad(\rho(\alpha))\gamma_*(u)\ra$. The
monodromy $N=M^{\nabla_{M\times S^1},\mu_{M\times S^1}}$ satisfies
$$N|_{X\times\{\alpha\}}=M_{\rho(\alpha)\gamma\rho(\alpha)^{-1}}.$$
Applying Corollary \ref{Mconstant} to $N$, we deduce our result.
\end{pf}

\section{The choice of the connection}
Recall that the map $a_1:H_1(G;\ZZ)\to H_1(X;\ZZ)$ is induced from the map
$a_1(x):G\ni g\mapsto gx\in X$, where $x\in X$ is an arbitrary point.
Through this section we will make the following topological assumption:
\begin{equation}
\Im a_1\cap\Tor H_1(X;\ZZ)=0. 
\label{condtop}
\end{equation}

Let $\alpha\in\Omega^2(X;\imag\RR)^G$ an invariant $2$-form representing
$-2\pi\imag c_1(L)$. Let $\mu\in\Omega^0(X;\imag\glie^*)^G$ satisfy
$d\mu(s)=\iota_{\fX(s)}\alpha$ for any $s\in\glie$, so that
$\alpha-\mu\in\Omega^2_G(X;\imag\RR)$ is a closed form in the Cartan complex,
and hence represents an equivariant cohomology class $[\alpha-\mu]\in
H^2_G(X;\imag\RR)$. 

\begin{lemma}
Suppose that $\frac{\imag}{2\pi}[\alpha-\mu]\in H^2_G(X;\ZZ)$. 
Then one can chose a connection $\nabla$ on $L$ whose 
curvature is $\alpha$ and such that for any morphism $\gamma:S^1\to G$
we have $M^{\nabla,\mu}_{\gamma}=1$.
More preciely, the set $\TT_{\alpha}^G(\mu)$ of gauge equivalence classes of 
connections satisfying this property is a torus of dimension 
$b_1(X)-\dim(\Im a_1\otimes_{\ZZ}\RR)$.
\label{onecanchose}
\end{lemma}
\begin{pf}
Let $\AAA_{\alpha}$ be the set of connections on $L$ whose curvature
is $\alpha$. Let $T\subset G$ be a maximal torus. 
By Lemma \ref{conjugats} it is
enough to consider $M_{\gamma}$ for $\gamma:S^1\to T$, since for any 
$\gamma:S^1\to G$ there exists $g\in G$ such that $g\gamma
g^{-1}(S^1)\subset T$. 

Let $\tlie=\Lie T$, and let $\Lambda=\Ker(\exp:\tlie\to T)$, so that
$T=\klie/\Lambda$. The morphisms $\gamma:S^1\to T$ are in 1--1
correspondence with elements of $\Lambda$.
For any $\gamma,\gamma'\in\Lambda$ we have 
\begin{equation}
M^{\nabla,\mu}_{\gamma'+\gamma}=
M^{\nabla,\mu}_{\gamma'}M^{\nabla,\mu}_{\gamma}.
\label{linearitat}
\end{equation}
To see this, observe the following.
Let $y\in L$, and let
$\nu(y;\cdot):\tlie\to L$ be the map defined as follows. For any $s\in\tlie$,
let $\nu^y_s:[0,1]\to L$ be the path such that $\nu^y_s(0)=y$ and
${\nu^y_s}'=\widetilde{\fX}^{\nabla,\mu}(s)(\nu^y_s)$. Then we set
$\nu(y;s)=\nu^y_s(1)$. With this definition, if $s\in\tlie$
and $v\in T_s\tlie\simeq\tlie$ (use the canonical isomorphism) then 
$D\nu(y;s)(v)=\widetilde{\fX}^{\nabla,\mu}(v)(s)$ (this is a consequence
of $[\widetilde{\fX}^{\nabla,\mu}(v),\widetilde{\fX}^{\nabla,\mu}(v')]=0$
for any $v,v'\in\tlie$). From this it follows that
\begin{equation}
\nu(\nu(y;s);s')=\nu(y;s+s'),
\label{esrep}
\end{equation}
which clearly implies (\ref{linearitat}).

Consider now the map $$c:\Lambda\to H^1(X;\RR)$$ which sends
$\gamma\in\Lambda$ to the homology class $[\gamma(S^1)]$ represented
by any orbit of the $S^1$ action on $X$ induced by $\gamma:S^1\to T$.
Let $\Lambda_0=\Ker c$. Using condition (\ref{condtop}), we deduce
from Lemma \ref{homognul} that for any
$\gamma\in\Lambda_0$ and any connection $\nabla\in\AAA_{\alpha}$ 
we have $M^{\nabla,\mu}_{\gamma}=1$. 
(Note that the map $\gamma^*:H^2_G(X;\ZZ)\to H^2_{S^1}(X;\ZZ)$ induced
by $\gamma$ lifts to the Cartan complex as
$\gamma^*(\alpha-\mu)=\alpha-\gamma^*(\mu)$.
Let now
$\Lambda_1=\Lambda_0^{\perp}$. This is a free abelian module. Let
$e_1,\dots,e_r\in\Lambda_1$ be a basis. By (\ref{linearitat}), if
a connection $\nabla\in\AAA_{\alpha}$
satisfies $M^{\nabla,\mu}_{e_j}=1$ for
any $1\leq j\leq r$, then $M^{\nabla,\mu}_{\gamma}=1$ for all
$\gamma\in\Lambda$. 

Finally, by gauge invariance of $M_{\gamma}^{\nabla,\mu}(x)$
(\ref{gaugeinvariant}), we can consider gauge classes of connections
on $L$ rather than connections. So let
$\TT_{\alpha}=\AAA_{\alpha}/\Map(X,S^1)$ be the gauge equivalence
classes of connections on $L$ with curvature $\alpha$. Picking a base
connection $\nabla\in\AAA_{\alpha}$, we can identify 
$$T_{\alpha}=\nabla+H^1(X;\RR)/H^1(X;\ZZ).$$
Furthermore, formula (\ref{addicio}) implies that if
$\eta\in\Omega^1(\imag\RR)$ satisfies $d\eta=0$, then 
$$M^{\nabla,\mu}_{e_j}=M^{\nabla+[\eta],\mu}_{e_j}-\la
[\eta],c(e_j)\ra,$$
where $[\eta]\in H^1(X;\RR)$ is the class represented by $\eta$.
On the other hand, 
the images by $c$ of $e_1,\dots,e_r$ are all linearly
independent. Hence,
$\la c(e_1),\dots,c(e_r)\ra$ is a space of dimension $r$, so 
the set of gauge equivalence classes of
connections $[\nabla]\in \TT_{\alpha}^G(\mu)$ 
such that $M^{\nabla,\mu}_{\gamma}=1$ for all $\gamma$ 
is the image under the quotient
$$H^1(X;\RR)\to H^1(X;\RR)/H^1(X;\ZZ)$$
of an affine subspace of codimension $r$. 
On the other hand, we have $r=\dim(\Im a_1\otimes\RR)$ so
$r\leq b_1(X)=\dim H^1(X;\RR)$, and hence this set is nonempty. More
precisely, the set $\TT_{\alpha}^G(\mu)\subset\TT_{\alpha}$ is a torus
of dimension $b_1(X)-r\geq 0$.
\end{pf}

\begin{lemma}
The infinitesimal lift $\widetilde{\fX}$ defined by $(\nabla,\mu)$
exponentiates to give a linear action of $G$ on $L$ if and only if for
any representation $\gamma:S^1\to G$ we have $M_{\gamma}=1$.
\label{sasufi}
\end{lemma}
\begin{pf}
Given $y\in L$ and $g\in G$, we define $gy\in L$ in the obvious way: let
$g=\exp(s)$, where $s\in\glie$, let $\nu^y_s:[0,1]\to L$ be the integral
curve of the vector field $\widetilde{\fX}^{\nabla,\mu}(s)$ with
initial value $\nu^y_s(0)=y$; then $gy:=\nu(y;s)=\nu^y_s(1)$. 

There are two things to check: that $gy$ is well defined and that the
resulting map $G\times L\to L$ is indeed an action of $G$ on $L$.
Observe first that both things are clear when $G=T$ is a torus (see
formula (\ref{esrep})). We now sketch how to deal with
the general case. Suppose that $s,s'\in\glie$ satisfy
$\exp(s)=\exp(s')$. We want to check that, for any $y\in L$,
$\nu^y_s(1)=\nu^y_{s'}(1)$. Now, it is easy to prove that there exists
some $s''\in\glie$ such that $\exp(s'')=\exp(s)=\exp(s')$ and such
that $[s,s'']=[s',s'']=0$. Then, applying (\ref{esrep}) to some torii
$T,T'$ such that $s,s''\in\Lie T$ and $s',s''\in\Lie T'$,
we deduce that $\nu(y;s)=\nu(y;s'')=\nu(y;s')$. This proves well
definedness. Finally, by Baker--Campbell--Haussdorf, $\nu$ satisfies
$\nu(\nu(y;s);s')=\nu(y;\log(\exp(s)\exp(s')))$ for $s,s'$ small
enough, and from this it follows easily that $\nu$ defines an action
of $G$ on $L$.
\end{pf}

\section{Proofs of the results}

We prove the theorems in two steps. First we assume that condition
(\ref{condtop}) is satisfied. Then we deduce the results in the general
case.

\subsection{Proofs of the theorems when (\ref{condtop}) holds}

\subsubsection{Proof of Theorem \ref{descripcio}}
Combine Lemma \ref{sasufi} with the ideas at the end of the proof of
Lemma \ref{onecanchose}.

\subsubsection{Proof of Theorem \ref{lifting}}

If the action of $G$ lifts to $L$, then the
first equivariant Chern class $c_1^G(L)$ of $L$ is an integral class 
and provides a lift of $c_1(L)$.
Now suppose that $c_1(L)=\iota^*(l)$, where $l\in
H^2_G(X;\ZZ)$. Take $\frac{\imag}{2\pi}(\alpha-\mu)\in\Omega^2_G(X)$ whose
cohomology class is equal in $H^2_G(X;\RR)$ to $l$. 
By Theorem \ref{descripcio}, there is some connection $\nabla$ on
$L$ which defines a lift of the action of $G$ to $L$. Let $L_{\alpha}$
denote the line bundle $L$ with the action of $G$.
Applying the Chern--Weil construction to equivariant bundles as
defined in \cite{BV} for the connection $\nabla$, we deduce that the form
$\frac{\imag}{2\pi}(\alpha-\mu)$ represents $c_1^G(L)$.
Now, since we have used de Rham theory, we have lost control of torsion, so 
that all we know in principle is that 
$$c_1^G(L)-l\in\Tor H^2_G(X;\ZZ).$$
To deduce that $c_1^G(L)=l$, we observe that the restriction of 
$\iota^*$ to $\Tor H^2_G(X;\ZZ)$ is an injection (indeed, 
$\Tor H^2_G(X;\ZZ)=\Ext (H_1(X_G;\ZZ),\ZZ)$,
$\Tor H^2(X;\ZZ)=\Ext (H_1(X;\ZZ),\ZZ)$ and, since $G$ is connected,
$\pi_1(BG)=0$, so the long exact sequence of homotopy groups for 
$X\to X_G\to BG$ tells us that $\pi_1(X)\to\pi_1(X_G)$ is exhaustive).
So, from $\iota^*(c_1^G(L)-l)=0$ we deduce that $c_1^G(L)=l$ in
$H^2_G(X;\ZZ)$.

To prove that $\iota^{-1}(c_1(L))$ classifies the lifts of the action to $L$
it is enough to check that if $G$ acts on $L$ and $c_1^G(L)=0$, then $L$ can be
$G$-equivariantly trivialised, i.e., there is an equivariant nowhere vanishing
section of $L$. So assume that $G$ acts on $L$ and $c_1^G(L)=0$. Take a
$G$-invariant connection $\nabla$ on $L$, let $\alpha=\nabla^2$ and let $\mu$ 
be the map given by Theorem \ref{tinflift}.
Now, by assumption $[\alpha-\mu]=0\in H^2_G(X;\imag\RR)$. Since the set of 
forms representing a fixed cohomology class is connected, we can join 
$\alpha-\mu$ to $0\in\Omega^2_G(X;\imag\RR)$ through a path 
$\gamma\subset\Omega^2_G(X;\imag\RR)$ all of
whose forms represent $0\in H^2_G(X;\imag\RR)$. 
Fix a trivialisation $L\simeq X\times\CC$.
It is easy to see, using 
the proof of Lemma \ref{onecanchose}, that $\gamma$ can be lifted continuously
to give $\forall t$ a connection $\nabla_t$ defining a lift to $L$ of the 
action with Chern--Weil form equal to $\gamma(t)$, in such a way that
$\nabla_1$ is the trivial connection on $L$. So we get a homotopy 
between the initial action of $G$ on $L$ and the trivial action defined from
a trivialisation $L\simeq X\times\CC$. Since $G$ is compact, this
implies that the initial action of $G$ on $L$ is trivial.

\subsection{Proof of the theorems in the general case}

Suppose that $T=\Im a_1\cap \Tor H_1(X;\ZZ)$ is nonzero.
Let $T'\subset H_1(X;\ZZ)$ be a complementary submodule of $T$, 
and let $G_T$ be the connected Lie group which fits in the exact
sequence
$$1\to T\to G_T\stackrel{q}{\to}G\to 1$$
with $q_*H_1(G_T;\ZZ)=T'$. The action of $G$ induces an action of 
$G_T$ on $X$, and we have a commutative diagram
$$\xymatrix{
H^2_G(X;\ZZ) \ar[rr]^{q}\ar[rd]_{\iota^*}
&& H^2_{G_T}(X;\ZZ) \ar[ld]^{\iota_T^*} \\
& H^2(X;\ZZ).}$$
Now, the action of $G_T$ clearly satisfies condition (\ref{condtop}), 
so we may apply the results obtained in the preceeding subsection and
get lifts of the action of $G_T$ to $L$, together with invariant 
connections. 

To prove Theorems \ref{lifting} and \ref{descripcio} for the action of
$G$ it is enough to check that, if $L_{\alpha}$ is a $G_T$ bundle 
isomorphic to $L$ (as bundles over $X$) such that
$$c_1^{G_T}(L_{\alpha})\in q^*H^2_G(X;\ZZ)$$
then the action of $G_T$ on $L_{\alpha}$ descends to an action of $G$,
or, equivalently, the action of $T\subset G_T$ on $L_{\alpha}$ is trivial
(note that, on the other hand, $q^*$ is injective).
This follows from the sequence of maps
$$H^2_G(X;\ZZ)\stackrel{q^*}{\to}
H^2_{G_T}(X;\ZZ)\stackrel{r^*}{\to}H^2(BT;\ZZ),$$
which is induced by the fibration $BT\to X_{G_T}\to X_G$, and 
consequently satisfies $r^*q^*=0$. The map $r^*$ is obtained from the 
$T$-equivariant inclusion $x_0\to X$ (where $x_0$ is any point).
And, since a representation $\rho:T\to\CC^*$ is trivial if and only
if $c_1^T(ET\times_{\rho}\CC)=0$, we deduce that if
$c_1^{G_T}(L_{\alpha})\in q^*H^2_G(X;\ZZ)$ then $r^*c_1^{G_T}(L_{\alpha})=0$
and hence $T$ acts trivially on $L_{\alpha}$.

\subsection{Proof of Corollary \ref{simplectic0}}

A theorem of Kirwan (see Proposition 5.8 in
\cite{Ki}) says that if $G$ acts in a
Hamiltonian fashion on $X$ then there is an isomorphism
$H^*_G(X;\QQ)\simeq H^*(X;\QQ)\otimes H^*(BG;\QQ)$.
In particular, this means that there exists an integer $d\geq 1$
such that if $a\in H^2(X;\ZZ)$ then $da\in\iota^* H^2_G(X;\ZZ)$.
So, given the line bundle $L$, there exists $l\in H^2_G(X;\ZZ)$ such
that $c_1(L^d)=\iota^*l$. Let now
$\frac{\imag}{2\pi}(\alpha'-\mu')\in\Omega^2_G(X)$ represent $l$ in
$H^2_G(X;\RR)$, and let $\nabla_d$ be a connection on $L^d$ whose
curvature is $\alpha'$. Let $\eta:=(\nabla^{\otimes d}-\nabla_d)^G$,
where $^G$ means the projection to the invariant subspace
$\Omega^1(X;\imag\RR)^G$ using the standard averaging trick:
$\zeta^G=\frac{1}{|G|}\int_{g\in G}g\zeta$. Then
$$\alpha-\mu:=(\alpha'-\mu')+d_{\glie}\eta\in\Omega^2_G(X;\imag\RR)$$
represents $-2\pi\imag l$, and the curvature of $\nabla$ is $\alpha$. 
On the other hand, $\Im a_1=0$, since any Hamiltonian action of 
$S^1$ on a compact manifold has fixed points, and hence the orbits
are contractible.
Consequently, by Theorem \ref{descripcio}, the lift
$\widetilde{\fX}^{\nabla^{\otimes d},\mu}$ exponentiates to an action
of $G$ on $L$ which leaves $\nabla$ fixed.

\subsection{Proof of Corollary \ref{simplectic}}
Let $\mu:X\to\glie^*$ be a moment map for the action of $G$ on $X$.
By Corollary \ref{simplectic0}, it suffices to take any
closed $\omega'-\mu'\in\Omega^2_G(X)$ 
near $\omega-\mu$ and representing a class 
$[\omega'-\mu']\in H^2_G(X;2\pi\QQ)$.

\end{document}